\documentclass{article}
\usepackage{spconf,amsmath,graphicx,bm,amssymb}
\usepackage{algorithm}
\usepackage{color}

\usepackage{graphicx,subfigure}

\usepackage[pdfstartview=FitH,
CJKbookmarks=true,
bookmarksnumbered=true,
bookmarksopen=true,
colorlinks,
pdfborder=001,
linkcolor=blue,
anchorcolor=blue,
citecolor=blue,
]{hyperref}
\hypersetup{hidelinks}

\usepackage{algorithmic}

\def\MatrixFont{\bf}
\def\VectorFont{\bf}

\newcommand{\mG}{{\MatrixFont G}}

\newcommand{\mI}{{\MatrixFont I}}

\newcommand{\mN}{{\MatrixFont N}}

\newcommand{\mQ}{{\MatrixFont Q}}
\newcommand{\mR}{{\MatrixFont R}}
\newcommand{\mS}{{\MatrixFont S}}

\newcommand{\mV}{{\MatrixFont V}}
\newcommand{\mW}{{\MatrixFont W}}
\newcommand{\mX}{{\MatrixFont X}}
\newcommand{\mY}{{\MatrixFont Y}}

\newcommand{\vh}{{\VectorFont h}}

\newcommand{\vn}{{\VectorFont n}}

\newcommand{\vs}{{\VectorFont s}}

\newcommand{\vw}{{\VectorFont w}}

\newcommand{\vy}{{\VectorFont y}}

\newtheorem{theorem}{Theorem}
\newtheorem{lemma}{Lemma}

\newtheorem{definition}{Definition}

\newcommand\tr{\ensuremath{{\rm Tr}}}

\title{Globally Optimal Beamforming Design for Integrated Sensing and Communication Systems}
\name{Zhiguo Wang$^{\star}$, Jiageng Wu$^+$,  Ya-Feng Liu$^{\S}$, and Fan Liu$^{\dag}$
}

\address{$^{\star}$College of Mathematics, Sichuan
University, Chengdu, China \\[2pt]
$^+$School of Mathematics, Jilin University, Changchun, China\\[2pt]
$^{\S}$LSEC, ICMSEC, AMSS, Chinese Academy of Sciences, Beijing, China\\[2pt]
    $^{\dag}$Southern University of Science and Technology, Shenzhen, China\\[2pt]
  Email:  wangzhiguo@scu.edu.cn, wujg22@mails.jlu.edu.cn, yafliu@lsec.cc.ac.cn, liuf6@sustech.edu.cn. }

\begin{document}

\ninept
\maketitle
%




%
%

\begin{abstract}
In this paper, we propose a multi-input multi-output beamforming transmit optimization model for joint radar sensing and multi-user communications, where the design of the beamformers is formulated as an optimization problem whose objective is  a weighted combination of the sum rate and the Cram\'{e}r-Rao bound,  subject to the transmit power budget constraint.  Obtaining a global solution for the formulated problem is a challenging task, because the sum rate maximization problem itself (even without considering the sensing metric) is known to be NP-hard.  In this paper, we propose an efficient global branch-and-bound algorithm for solving the formulated problem based on the McCormick envelope relaxation and the semidefinite relaxation technique. The proposed algorithm is guaranteed to find the global solution for the considered problem, and thus serves as an important benchmark for performance evaluation of the existing local or suboptimal algorithms for solving the same problem.

\end{abstract}
\vspace{-0.1cm}
\section{Introduction}\vspace{-0.2cm}
\label{sec:intro}
The integrated sensing and communication (ISAC) technology has gained significant attention and recognition from both academia and industry as a pivotal enabler \cite{ITU-R-WP5D2023,rahman2019framework,tan2018exploiting,ma2020joint,xiong2023fundamental,song2023intelligent}. Numerous studies have investigated the signaling strategy in multi-antenna ISAC systems, with a  particular focus on the joint beamforming optimization \cite{liu2022integrated,he2023full}. The recent work \cite{liu2021cramer} formulated a  Cram\'{e}r-Rao bound (CRB)-minimization problem  under the  signal-to-interference-plus-noise ratio (SINR) constraints \cite{liu2022joint}, which aim to guarantee a minimum level of communication quality of service at all users \cite{liu2020joint}.    Since the sum rate \cite{shi2011iteratively,shen2018fractional} is a more fundamental metric that comprehensively characterizes the overall performance in multi-user scenarios, this paper is interested in optimizing both the communication performance, measured by the sum rate for different users, and
the target estimation performance, measured by the CRB for unbiased estimators.
 These optimization objectives are subject to the  total transmit power budget constraint. In contrast to individual SINR constraints for each user in \cite{liu2020joint}, our design incorporates sum rate maximization (SRM) to ensure the superior network throughput performance.

However, the task of maximizing the sum rate poses greater difficulty compared to the fixed SINR target case.  Notably, the SRM problem has been proven to be NP-hard in \cite{chiang2007power,luo2008dynamic}. There are generally two popular methods for achieving the global solution to the SRM problem \cite{weeraddana2011weighted,matthiesen2020mixed}, where one is the outer polyblock approximation (PA) algorithm \cite{qian2009mapel}, and the other is the branch and bound (B\&B) method \cite{tuy2005monotonic}. In general, the key to the success of the PA algorithm depends on the property that  the objective function of SRM is monotonically increasing in SINRs \cite{liu2012achieving}, but this monotonicity property  does not hold for our interested problem, since the objective function also has a CRB term for radar sensing. As described in \cite{tuy2005monotonic}, the numerical efficiency as well as the convergence speed of the B\&B algorithm heavily relies on  the quality of upper and lower bounds for the optimal value \cite{weeraddana2011weighted,lu2017efficient,lu2020enhanced,schobel2010theoretical}.
 Recently,  the work in \cite{zhu2023information} utilized a combination of the weighted mean square error minimization (WMMSE) and the SDR techniques to achieve a Karush-Kuhn-Tucker point for the problem. But this algorithm is not guaranteed to find the globally optimal solution of the problem.

 In this paper, we consider the problem of optimizing the weighted combination of the sum rate of all communication users and the CRB of the sensing target under the total transmit power constraint. We propose a computationally  efficient global algorithm for solving the formulated nonconvex problem, which lies in the celebrated B\&B framework while is based on the McCormick envelope relaxation and the  semidefinite relaxation (SDR) technique. The key feature of the proposed algorithm is its global optimality guarantee. Indeed,  the computation of a global solution holds significant importance, as the resultant global solution plays a crucial role in evaluating the inherent performance limits of the associated wireless communication system. Moreover, global optimization algorithms serve as vital benchmarks for assessing the performance of existing heuristic or local algorithms designed for the same problem.
\vspace{-0.1cm}
\section{Problem Formulation}\label{sec_system}\vspace{-0.2cm}
\vspace{-0.0cm}
Consider a  multi-input multi-output (MIMO) ISAC base station (BS) equipped with $N_t$ transmit
antennas and $N_r$ receive antennas, which serves $K$
downlink single-antenna users while detecting an extended target. In this paper, we assume $K\leq N_t<N_r$ \cite{liu2021cramer} in order to guarantee the feasibility of the beamforming design problem while avoiding the information loss of the sensed target.

Let $\mX\in \mathbb{C}^{N_t\times L}$ be the transmitted baseband signal, which is  the sum of linear precoded radar waveforms and communication symbols, given by
\vspace{-0.3cm}
\begin{align}\label{eqn:X}
\mX=\sum_{k=1}^{K}\vw_k\vs_k+\mW_A\mS_A,
\end{align}
where $\vs_k\in\mathbb{C}^{1 \times L}$ is the data symbol for the $k$-th communication user and $\mS_A\in\mathbb{C}^{N_t\times L}$ is the sensing signal, which are precoded by
the communication beamformer $\vw_k\in\mathbb{C}^{N_t\times 1}$ and the auxiliary beamforming matrix $\mW_A\in\mathbb{C}^{N_t\times N_t}$, respectively.
Assume that the data streams $\tilde{\mS}=[\vs_1^H,\ldots,\vs_K^H,\mS_A^H]^H$ are asymptotically orthogonal \cite{liu2021cramer} to each other for sufficiently large $L$, i.e., $\frac{1}{L}\tilde{\mS}\tilde{\mS}^H\approx\mI_{K+N_t}$.

For multi-user communications, by transmitting $\mX$ to $K$ users, the received signal $\vy_k$ of user $k$ is given as
\begin{align}\label{eqn:rsm}
\vy_k = \vh_k^H\mX+\vn_C,
\end{align}
where $\vh_k\in \mathbb{C}^{ N_t\times 1}$ is the communication  channel between the BS and user $k$, which is assumed to be known to the BS; $\vn_C$ is an additive white Gaussian noise (AWGN) vector with the variance of
each entry being $\sigma_C^2$.
 Then the SINR at the $k$-th communication user can be expressed as
 \vspace{-1mm}
 \begin{align}\label{eqn:SINR}
  \tilde{\gamma}_k=\frac{\left|\vh_{k}^{H} \vw_{k}\right|^{2}}{\sum_{i=1, i \neq k}^{K}\left|\vh_{k}^{H} \vw_{i}\right|^{2}+\left\|\vh_{k}^{H} \mW_A\right\|^{2}+\sigma_{C}^{2}},~k\in[K].
 \end{align}
 where $[K]$ denotes the set $\{1,2,\ldots,K\}$.
One of the most important criteria for multiuser beamforming is the overall system throughput,
\vspace{-1mm}
\begin{align}\label{sum_rate}
\sum_{k=1}^{K}\log(1+\tilde{\gamma}_k).
\end{align}

By transmitting $\mX$ to sense the target, the reflected echo signal at the BS is given by
\vspace{-1mm}
\begin{align}\label{eqn: YR}
\mY_s=\mG\mX+\mN_s,
\end{align}
where $\mG\in \mathbb{C}^{N_r\times N_t}$ denotes the target response matrix; $\mN_s\in \mathbb{C}^{N_r\times L}$ is an AWGN matrix, with the variance of each entry being  $\sigma_s^2$.
For the purpose of target sensing, we focus on estimating the response
matrix $\mG$. In order to improve the estimation performance, the CRB of estimating the response matrix is given by
\vspace{-1mm}
\begin{align}\label{eqn:CRB}
\textmd{CRB}(\mG) = \frac{\sigma_s^2N_r}{L}\tr(\mR_X^{-1}),
\end{align}
where $\mR_X =\frac{1}{L}\mX\mX^H= \sum_{k=1}^{K}\vw_k\vw_k^H+\mW_A\mW_A^H$
is the sample covariance matrix of $\mX$.

Based on the sum rate expression in \eqref{sum_rate} and the CRB  expression in \eqref{eqn:CRB}, the joint communication and sensing beamforming design problem can be formulated as
\vspace{-2mm}
\begin{subequations}\label{beamform}
\begin{align}
\label{beamform_obj}&\min_{\substack{\{\vw_k\}_{k=1}^K,\mW_A}}~ -\sum_{k=1}^K\log(1+\tilde{\gamma}_k)+ \rho\tr\left(\mR_X^{-1}\right) \\
\label{beamform_cons2} &\qquad~\textmd{s.t.} ~ ~~~~ \tr\left(\sum_{k=1}^{K}\vw_k\vw_k^H+\mW_A\mW_A^H\right)\leq P_T,
\end{align}
\end{subequations}
where $\rho$ is a parameter to trade-off the sum rate and  the CRB; $P_T$ is the total transmit power budget
of the BS.
It is worth highlighting that, although we formulate the problem as in \eqref{beamform}, the proposed algorithm in this paper can also be used to solve the other formulations of the joint communication and sensing beamforming design problem such as the SRM problem subject to the sensing CRB constraint on the target and the total power budget constraint of the BS.
\vspace{-2mm}
\section{Proposed Global Algorithm}
In this section, we propose a global optimization algorithm for solving problem \eqref{beamform}. The proposed algorithm lies in the B\&B framework and the lower bound in our algorithm is based on the SDR technique and the McCormick envelope relaxation.
\vspace{-2mm}
\subsection{SDR and McCormick Envelope Relaxation of Problem \eqref{beamform}}
Introducing some auxiliary variables $\left\{\Gamma_k\right\}_{k=1}^K$,  we can reformulate problem \eqref{beamform} as
\vspace{-0.3cm}
 \begin{subequations}\label{beamform_multi}
\begin{align}
\label{beamform_multi_obj}&\min_{\substack{\{\Gamma_k\}_{k=1}^K,\\
\{\vw_k\}_{k=1}^K,\mW_A}}~ -\sum_{k=1}^K\log(1+\Gamma_k)+ \rho \tr\left(\mR_X^{-1}\right) \\
\label{beamform_multi_cons2} &\qquad\textmd{s.t.}  \qquad~~ \tr\left(\sum_{k=1}^{K}\vw_k\vw_k^H+\mW_A\mW_A^H\right)\leq P_T,\\
 \label{beamform_multi_cons3}&\qquad\qquad\qquad\tilde{\gamma}_k\geq \Gamma_k,~ k\in [K].
\end{align}
\end{subequations}
Let
$
\mW_k=\vw_k\vw_k^H$ for all $k\in[K]$ and $\mW_{K+1}=\mW_A\mW_A^H
$. Then we get $\mW_k\succeq 0$ and
 $\textmd{rank}(\mW_k)=1$ for all $k\in [K].$
According to the definitions of $\tilde{\gamma}_k$ in \eqref{eqn:SINR}, the $k$-th SINR constraint \eqref{beamform_multi_cons3} can be rewritten as
\begin{align}\label{eqn:QW}
\tr(\mQ_k\mW_k)-\Gamma_k\sum_{i\neq k}^{K+1}\tr(\mQ_k\mW_i)\geq \Gamma_k\sigma_C^2,~ k\in [K],
\end{align}
where $\mQ_k=\vh_k\vh_k^H$.  Then problem \eqref{beamform_multi} can be rewritten as
 \begin{subequations}\label{beamform_matrix}
\begin{align}
\label{beamform_matrix_obj}\min_{\substack{\{\Gamma_k\}_{k=1}^{K},\{\mW_k\}_{k=1}^{K+1}}}&~ -\sum_{k=1}^K\log(1+\Gamma_k)+ \rho \tr\left(\mR_X^{-1}\right) \\
\label{beamform_matrix_cons2} \textmd{s.t.}~~~\qquad&  ~ \sum_{k=1}^{K+1}\tr\left(\mW_k\right)\leq P_T,~ \eqref{eqn:QW},\\
  \label{beamform_matrix_cons4}&~\textmd{rank}(\mW_k)=1, \mW_k\succeq 0, ~k\in [K].
\end{align}
\end{subequations}
By dropping the rank constraints in \eqref{beamform_matrix_cons4}, problem \eqref{beamform_matrix} is relaxed into the following optimization problem
 \begin{subequations}\label{beamform_matrix1}
\begin{align}
\label{beamform_matrix_obj1}&\min_{\substack{\{\Gamma_k\}_{k=1}^{K},\{\mW_k\}_{k=1}^{K+1}}} -\sum_{k=1}^K\log(1+\Gamma_k)+ \rho \tr\left(\mR_X^{-1}\right) \\
\nonumber &\qquad~~~~~\textmd{s.t.}  \sum_{k=1}^{K+1}\tr\left(\mW_k\right)\leq P_T,~\eqref{eqn:QW},~\mW_k\succeq 0,  ~k\in [K].
\end{align}
\end{subequations}
Note that the above optimization problem \eqref{beamform_matrix1} is not convex
since there exists a bilinear term in \eqref{eqn:QW}. In the following, we shall use the McCormick envelope relaxation technique to deal with this nonconvexity. Before doing that, below we first consider how to construct a rank-one solution from the solutions of problem \eqref{beamform_matrix1} under the assumption that the optimal solution of problem \eqref{beamform_matrix1} can be obtained.

From the definition of $\mR_X$, problem \eqref{beamform_matrix1} can be rewritten as
\vspace{-2mm}
 \begin{subequations}\label{beamform_multi1}
\begin{align}
\label{beamform_multi_obj1}&\min_{\substack{\mR_X,\{\Gamma_k\}_{k=1}^K,\\
\{\mW_k\}_{k=1}^K}} -\sum_{k=1}^K\log(1+\Gamma_k)+ \rho \tr\left(\mR_X^{-1}\right) \\[-3mm]
\label{beamform_multi_cons12} &\qquad~\textmd{s.t.}  \qquad \tr(\mR_X)\leq P_T, \mR_X\succeq \sum_{k=1}^{K}\mW_k,
\\
\label{beamform_multi_cons13}&\qquad\qquad \qquad~\eqref{eqn:QW},~ \mW_k\succeq 0, k\in[K].
\end{align}
\end{subequations}
 There are a lot works \cite{godara2018handbook,pataki1998rank,ma2017unraveling} that study the tightness of SDRs in the context of beamformer design. It encourages us to believe that the relaxation used in \eqref{beamform_multi1} from \eqref{beamform_matrix} is tight.
\begin{theorem}\label{theorem_rank}
Given an optimal solution $\bar{\mR}_X, \{\bar{\Gamma}_k\}_{k=1}^K, \{\bar{\mW}_k\}_{k=1}^K$ of problem \eqref{beamform_multi1}, the following
$\tilde{\mR}_X, \{\tilde{\Gamma}_k\}_{k=1}^K, \{\tilde{\mW}_k\}_{k=1}^K$ is also an optimal solution:
\begin{align}
\tilde{\mR}_X= \bar{\mR}_X, \tilde{\Gamma}_k = \bar{\Gamma}_k,  \tilde{\mW}_k=\frac{\bar{\mW}_k\mQ_k\bar{\mW}_k^H}{\tr(\mQ_k\bar{\mW}_k)},~k\in[K].
\end{align}
Moreover,  $\textmd{rank}(\tilde{\mW}_k)=1$,  for all $k\in [K]$.
\end{theorem}
We can use the results in Theorem \ref{theorem_rank} to find a rank-one optimal solution $\tilde{\mR}_X, \{\tilde{\Gamma}_k\}_{k=1}^K, \{\tilde{\mW}_k\}_{k=1}^K$ if the optimal solution $\bar{\mR}_X, \{\bar{\Gamma}_k\}_{k=1}^K,\{\bar{\mW}_k\}_{k=1}^K$ of problem \eqref{beamform_multi1} is obtained.
In addition, the optimal beamformer $\vw_k$ for the original problem \eqref{beamform} is straightforwardly expressed as
\vspace{-3mm}
\begin{align*}
\vw_k=\left(\vh_k^H\bar{\mW}_k\vh_k\right)^{-1/2}\mW_k\vh_k, ~k\in [K]
\end{align*}
and the beamformer $\mW_A$ is calculated by the Cholesky decomposition
$
\mW_A\mW_A^H=\tilde{\mR}_X-\sum_{k=1}^{K}\tilde{\mW}_k.
$
Therefore, we only need to  focus on solving problem \eqref{beamform_multi1} in order to solve problem \eqref{beamform_multi}. Finally, Theorem \ref{theorem_rank} can be proved by using the similar argument as in \cite{liu2020joint}.

Now we deal with the nonconvexity coming from the bilinear term in constraint \eqref{eqn:QW}. To do so,
let us introduce some auxiliary variables $s_k$, $k\in[K]$. Then the $k$-th  nonconvex constraint \eqref{beamform_multi_cons13} can be rewritten as
\vspace{-2mm}
\begin{align}
\label{eqn_sk} \tr(\mQ_k\mW_k)-s_k\geq \Gamma_k\sigma_C^2, ~
 s_k = \Gamma_k\tr(\mQ_k(\mR_X-\mW_k)).
\end{align}
One can observe that  there is still a bilinear function in \eqref{eqn_sk}.
Next, we develop a  convex relaxation for \eqref{eqn_sk} based on the McCormick envelopes \cite{mitsos2009mccormick}.
\begin{lemma}\label{lemma_ebve}
Assuming that $\ell_k\leq\Gamma_k\leq u_k$ and $0\leq\tr(\mQ_k(\mR_X-\mW_k))\leq b_k$, $k\in[K]$, then the McCormick envelopes for the bilinear constraint \eqref{eqn_sk} is
\vspace{-2mm}
\begin{subequations}\label{enve_lopes}
\begin{align}
\label{enve_lopes1} s_k \geq& ~\ell_k\tr(\mQ_k(\mR_X-\mW_k)),\\
\label{enve_lopes2} s_k \geq&~ u_k\tr(\mQ_k(\mR_X-\mW_k))+(\Gamma_k-u_k)b_k,\\
\label{enve_lopes3} s_k \leq&~ u_k\tr(\mQ_k(\mR_X-\mW_k)),\\
\label{enve_lopes4} s_k \leq&~ (\Gamma_k-\ell_k)b_k+\ell_k\tr(\mQ_k(\mR_X-\mW_k)),
\end{align}
\end{subequations}
all of which are linear constraints with respect to $\mR_X, \{\Gamma_k\}_{k=1}^K$, and $\{\mW_k\}_{k=1}^K$.
\end{lemma}
From the power constraint $\tr(\mR_X)\leq P_T$, it is simple to derive that
\vspace{-5mm}
\begin{align*}
0\leq\Gamma_k\leq \frac{P_T\|\vh_k\|^2}{\sigma_C^2}, \tr(\mQ_k(\mR_X-\mW_k))\leq P_T\|\vh_k\|^2,~ k\in[K].
\end{align*}
The above lower and upper bounds on the SINR and interference terms provide desirable bounds in Lemma \ref{lemma_ebve}.
Setting $b_k=P_T\|\vh_k\|^2$ in the above Lemma \ref{lemma_ebve}, we immediately obtain the following convex McCormick envelope based relaxation (MER) of problem \eqref{beamform_multi1}:
\vspace{-3mm}
 \begin{subequations}\label{MER_multi1}
\begin{align}
\label{MER_multi_obj1}&\min_{\substack{\mR_X,\{\Gamma_k\}_{k=1}^K,\\
\{\mW_k\}_{k=1}^K,\{s_k\}_{k=1}^K}} -\sum_{k=1}^K\log(1+\Gamma_k)+ \rho \tr\left(\mR_X^{-1}\right)\\[-2mm]
\label{MER_multi_cons12} &\qquad~~~~\textmd{s.t.} ~ \tr(\mR_X)\leq P_T, \mR_X\succeq \sum_{k=1}^{K}\mW_k, \mW_k\succeq 0,\\
\label{MER_multi_cons13}&\qquad\qquad~\tr(\mQ_k\mW_k)-s_k\geq \Gamma_k\sigma_C^2, \ell_k\leq\Gamma_k\leq u_k,\\
\label{MER_multi_cons14}& \qquad\qquad~0\leq\tr(\mQ_k(\mR_X-\mW_k))\leq b_k,~\eqref{enve_lopes},~k\in[K],
\end{align}
\end{subequations}
which is a convex problem that can be solved via CVX \cite{grant2014cvx}.
\vspace{-3mm}
\subsection{Proposed B\&B Algorithm}
Now we are ready to present our proposed B\&B algorithm for globally solving problem \eqref{beamform_multi1}. The basic idea of the proposed
algorithm is to relax the original problem \eqref{beamform_multi1} (with bilinear constraints) to MER \eqref{MER_multi1} and gradually tighten the relaxation by reducing the width of the associated intervals $[\ell_k,u_k]$ for the $k$-th communication user.

For ease of presentation, we first introduce some notations. Let MER($\mathcal{Q}$) denote as the corresponding MER problem defined over the rectangle set $\mathcal{Q}:=\prod_{k=1}^K[\ell_k,u_k]$; let $L$ be the optimal value of MER($\mathcal{Q}$) and $\mV:=\left[\{\Gamma_k\}_{k=1}^K,\{s_k\}_{k=1}^K,\{\mW_k\}_{k=1}^K,\mR_X\right]$ be its
optimal solution; let  $\mathcal{P}$ denote the constructed problem list and $\{\mathcal{Q},\mV, L\}$ denote a problem instance from
the list $\mathcal{P}$; let $U^t$ denote the upper bound at the $t$-th iteration; let $\mV^*$ denote the best known feasible solution and $U^*$ denote the
objective value of problem \eqref{beamform_multi1} at $\mV^*$.
Then we can present the following key components of the proposed B\&B algorithm.

 \emph{Initialization}: We initialize all intervals $[\ell_k^0,u_k^0]$ for all $k\in[K]$ to be $\left[0, \frac{P_T\|\vh_k\|^2}{\sigma_C^2}\right]$ for the MER problem, and $\mathcal{Q}^0:= \prod_{k=1}^K[\ell_k^0,u_k^0]$. In this
case, we use the CVX  to solve the corresponding problem \eqref{MER_multi1}  and obtain its optimal solution $\mV^0$ and its optimal value $L^0$.

 \emph{Termination}: Let $\{\mathcal{Q}^t,\mV^t, L^t\}$ denote the problem instance that has the least lower bound in the problem list $\mathcal{P}$. If
\begin{align}\label{eqn:term}
U^t-L^t\leq \epsilon,
\end{align}
where $\epsilon$ is the given error tolerance, we stop the algorithm; otherwise we branch one interval in $\mathcal{Q}^t$ as specified below in \eqref{eqn:k}.

 \emph{Branch}: Suppose that the stopping criterion in \eqref{eqn:term} is not satisfied, we select one interval that leads to the largest relaxation gap to be branched to smaller sub-intervals. Let $\mV^t$ be the optimal solution of problem MER($\mathcal{Q}^t$). Since problem MER($\mathcal{Q}^t$) is a relaxation of problem \eqref{beamform_multi1}, its solution might not satisfy the constraint \eqref{beamform_multi_cons13}. Fortunately, we can construct a feasible solution to problem \eqref{beamform_multi1} based on the solution of MER($\mathcal{Q}^t$) as follows:
    \vspace{-2mm}
     \begin{subequations}\label{eqn:feas}
    \begin{align}
    &\hat{\mR}_X^t= \mR_X^t,~\hat{\mW}_k^t=\mW_k^t,~k\in[K],\\
&\hat{\Gamma}_k^t~~ =  \frac{\Gamma_k^t\sigma_C^2+\ell_k^t\tr(\mQ_k(\mR_X^t-\mW_k^t))}{\sigma_C^2+\tr(\mQ_k(\mR_X^t-\mW_k^t))}, ~k\in[K].
\end{align}
\end{subequations}
 The  constructed solution in \eqref{eqn:feas} plays a central role in improving the upper bound in the algorithm and in selecting the user that leads to the largest relaxation gap to be branched.

In particular, we use the following rule to  select the user that has the largest relative relaxation gap:
\vspace{-2mm}
\begin{align}\label{eqn:k}
 k^* &=\mathop{\arg\max}\limits_{k\in \{1,2,\ldots,K\}} \left\{ \frac{\Gamma_k^t-\hat{\Gamma}_k^t}{1+\hat{\Gamma}_k^t}\right\}.
\end{align}
It is clear that the numerator in \eqref{eqn:k} is the gap between the predicted SINR (by the relaxation) and the practically achievable SINR of user $k$ and hence the quantity in \eqref{eqn:k} measures the relative relaxation gap between the predicted and achievable SINRs of user $k.$

 Next,
 we  partition $\mathcal{Q}^t$ into two sets (denoted as $\mathcal{Q}_1^t$ and $\mathcal{Q}_2^t$) by partitioning its $k^*$-th interval into two equal intervals and keep all the others being unchanged. Then we solve the MER problems over the newly obtained two small sets, which are called children problems. Obviously, the two children problems obtained from partitioning $\mathcal{Q}^t$ are tighter than the one defined over the original set $\mathcal{Q}^t$. In this way, the B\&B process gradually tightens the relaxations and is able to find a (nearly) global solution satisfying the condition in \eqref{eqn:term}.
    When $\mathcal{Q}^t$ has been branched into two sets, the problem instance defined over $\mathcal{Q}^t$ will be deleted from the problem list $\mathcal{P}$, and the two corresponding children problems will be added into $\mathcal{P}$ if their optimal objective values are less than or equal to the current upper bound.

\emph{Lower Bound}: For any problem instance $\{\mathcal{Q},\mV, L\}$, $L$ is the lower bound of the optimal value of the original nonconvex problem \eqref{beamform_multi1} defined over $\mathcal{Q}$.
Therefore, the smallest one among all bounds is a lower bound of the optimal value of the original problem \eqref{beamform_multi1}.  At the $t$-th iteration, we choose a problem instance from $\mathcal{P}$, denoted as $\{\mathcal{Q}^t,\mV^t,L^t\}$, such that the bound $L^t$ is the smallest one in $\mathcal{P}$.

\emph{Upper Bound}: An upper bound is obtained from the best known feasible solution of \eqref{beamform_multi1}.
Following \eqref{eqn:feas}, $\hat{\mR}_X^t$, $\{\hat{\Gamma}_k^t\}_{k=1}^K$, $\{\hat{\mW}_k^t\}_{k=1}^K$ is a feasible solution for problem \eqref{beamform_multi1} and hence
\vspace{-3mm}
\begin{align}\label{eqn:upper_obj}
\hat{U}^t:=-\sum_{k=1}^K\log(1+\hat{\Gamma}_k^t)+ \rho \tr((\hat{\mR}_X^t)^{-1})
\end{align}
is an upper bound of the original problem. In our proposed algorithm, the upper
bound $U^t$ is the best objective values at all of the known feasible
solutions at the $t$-th iteration.

The pseudo-codes of our proposed algorithm are given in Algorithm \ref{alg_BB}, which is a careful combination of all of the above key components.

 To the best of our knowledge, our proposed algorithm is the first global algorithm for solving problem \eqref{beamform_multi1}. Before presenting the main theoretical result,
\begin{algorithm}[t]
\begin{algorithmic}
\STATE {\bf Initialization} Give an error tolerance $\epsilon>0$. Solve problem MER$(\mathcal{Q}^0)$ \eqref{MER_multi1} at $t=0$ to obtain  $\{\mathcal{Q}^0,\mV^0,L^0\}$, and add it into the problem list $\mathcal{P}$.  Compute a feasible point $\left\{\hat{\mR}_X^0, \{\hat{\Gamma}_k^0\}_{k=1}^K,\{\hat{\mW}_k^0\}_{k=1}^K \right\}$ and an upper bound
$U^0=\hat{U}^0$ by \eqref{eqn:feas} and \eqref{eqn:upper_obj}, respectively.
\FOR{ $t=1,2,\ldots$}
\IF{$U^t-L^t<\epsilon$}
\STATE terminate the algorithm and return $U^*=U^t$ and $\mV^*=\mV^t$.
\ENDIF
\STATE  Choose the problem that has the lowest bound from $\mathcal{P}$ and
delete it from $\mathcal{P}$.
\STATE Obtain $k^*$ by \eqref{eqn:k} and compute $z_{k^*}^t=\frac{1}{2}(\ell_{k^*}^t+u_{k^*}^t)$.
\STATE Branch $\mathcal{Q}^t$ into two sets $\mathcal{Q}_1^t=\{\Gamma_k\in \mathcal{Q}^t\mid\Gamma_{k^{*}}\leq z_{k^*}\}$
and $\mathcal{Q}_2^t=\{\Gamma_k\in \mathcal{Q}^t\mid\Gamma_{k^{*}}\geq z_{k^*}\}$
\FOR{ $j=1,2$}
\STATE Solve problem MER \eqref{MER_multi1}  with set $\mathcal{Q}_j^t$  to obtain  $\mV_j^t$ and  $L_j^t$.
\STATE Compute  $\{\hat{\mR}_{X,j}^t, \{\hat{\mW}_{k,j}^t\}_{k=1}^K,\{\hat{\Gamma}_{k,j}^t\}_{k=1}^K\}$ and  $\hat{U}^t$ in \eqref{eqn:feas} and \eqref{eqn:upper_obj}, respectively.
\IF{$L_j^t\leq U^t$}
\STATE add $\{\mathcal{Q}^t,\mV_{j}^t,L_j^t\}$ into the problem list $\mathcal{P}$.
\ENDIF
\IF{$U^t\geq \hat{U}^t$}
\STATE set $U^t=\hat{U}^t$, $ \hat{\mR}_{X}^*=\hat{\mR}_{X,j}^t,\hat{\mW}_{k}^*=\hat{\mW}_{k,j}^t,\Gamma_k^*=\hat{\Gamma}_{k,j}^t$.
\ENDIF
\ENDFOR
\ENDFOR
\end{algorithmic}
\caption{Proposed B\&B Algorithm for Solving Problem \eqref{beamform_multi1}}
\label{alg_BB}
\end{algorithm}
let us first formally define the
$\epsilon$-optimal solution of problem \eqref{beamform_multi1}.
\vspace{-2mm}
\begin{definition}
 Given any $\epsilon>0$, a feasible point $\mR_X$, $\{\Gamma_k\}_{k=1}^K$, $\{\mW_k\}_{k=1}^K$ is called an $\epsilon$-optimal solution of problem \eqref{beamform_multi1} if it satisfies
 \vspace{-4mm}
\begin{align}\label{eop_solution}
\Phi\left(\mR_X,\{\Gamma_k\}_{k=1}^K,\{\mW_k\}_{k=1}^K\right)-\nu^*\leq \epsilon,
\end{align}
where
$
\Phi\left(\mR_X,\{\Gamma_k\}_{k=1}^K,\{\mW_k\}_{k=1}^K\right):=-\sum_{k=1}^K\log(1+\Gamma_k)+ \rho \tr\left(\mR_X^{-1}\right)
$
is the objective value of problem \eqref{beamform_multi1}  at the point $\mR_X$, $\{\Gamma_k\}_{k=1}^K$, $\{\mW_k\}_{k=1}^K$  and $\nu^*$ is the
optimal value of problem \eqref{beamform_multi1}.
\end{definition}
\vspace{-3mm}
\begin{theorem}\label{theorem_global}
For any given $\epsilon>0$ and any given instance of problem \eqref{beamform_multi1} with $K$ users,  Algorithm \ref{alg_BB} will terminate and return an $\epsilon$-optimal solution defined in \eqref{eop_solution} within
$
\left(\frac{\Gamma_{\max}}{\delta}\right)^{K}+1
$
iterations, where $\Gamma_{\max}= \mathop{\max}\limits_{k\in [K]}\frac{P_T\|\vh_k\|^2}{\sigma_C^2}$ and $\delta=\frac{1}{2}\left(\exp\left(\frac{\epsilon}{K}\right)-1\right).$
\end{theorem}
 Theorem \ref{theorem_global} shows that
the total number of iterations for our proposed  algorithm to return an $\epsilon$-optimal solution grows exponentially fast with the total number of users $K.$ The iteration complexity of the proposed algorithm seems to be prohibitively high. However, our simulation results in the next section show that its practical iteration complexity is actually significantly less than the worst-case bound in Theorem \ref{theorem_global}, thanks to the effective lower bound provided by the relaxation problem \eqref{MER_multi1} based on the SDR technique and the McCormick envelope relaxation.
\vspace{-0.3cm}
\section{Numerical Results}\label{sec_num}
 Consider a MIMO BS that is equipped with $N_t=10$ antennas serving $K=5$ users. Set $\sigma_C^2=0.1$ and assume that the  channel $\vh_k$ for the $k$-th user follows the standard complex  Gaussian distribution.

\begin{figure}[t]
\begin{center}
\includegraphics[width=0.75\linewidth]{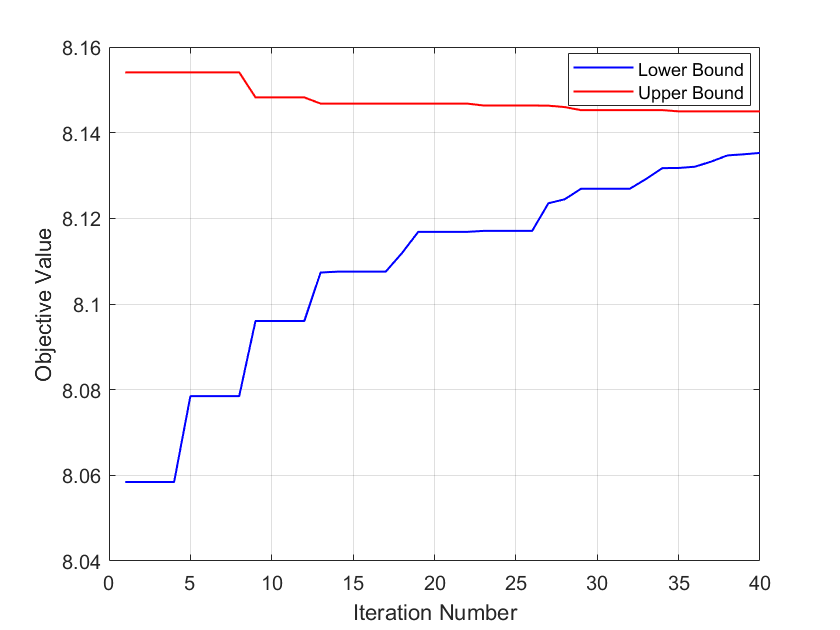}
\vspace{-0.4cm}
\caption{Lower and upper bounds of the proposed B\&B algorithm. }
\label{lower_upper_new}
\end{center}
\vspace{-0.7cm}
\end{figure}

 The convergence behavior of our proposed algorithm is shown in Fig. \ref{lower_upper_new}. In this simulation, we set
the error tolerance $\epsilon=0.01$, the power budget $P_T=10$ dB, and the weight parameter $\rho=1$ in Algorithm \ref{alg_BB}.  The result is averaged over 10 Monte Carlo runs. It can be seen from Fig. \ref{lower_upper_new} that the proposed B\&B algorithm quickly converges to the global solution. In particular, the total number of iterations for the proposed algorithm to find the $\epsilon$-optimal solution is $40$ (and the CPU time is $139$ seconds).

\begin{figure}[t]
\begin{center}
\includegraphics[width=0.75\linewidth]{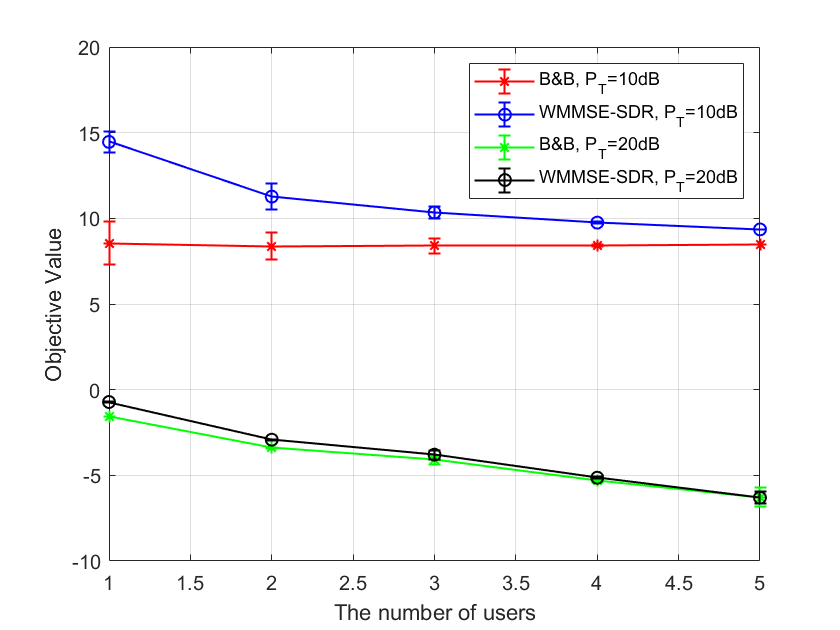}
\vspace{-0.4cm}
\caption{Comparison of the proposed B\&B algorithm with the WMMSE-SDR algorithm \cite{zhu2023information} in terms of the objective value.}
\label{compare_local}
\end{center}
\vspace{-0.8cm}
\end{figure}

Furthermore, we evaluate the performance of an existing state-of-the-art algorithm \cite{zhu2023information} for solving the original nonconvex problem \eqref{beamform},
which is called WMMSE-SDR. However, the WMMSE-SDR algorithm is only guaranteed to find a stationary point of problem \eqref{beamform}. In Fig. \ref{compare_local}, we show the impact of the communication user number on the objective value, with the power budget being set as 10 dB and 20 dB, respectively.  It can be observed from Fig. \ref{compare_local} that
the gap of the objective values at the solutions found by  WMMSE-SDR and the proposed B\&B algorithm becomes larger when the power increases. The results demonstrate the performance gain of the proposed global optimization algorithm over that of the local optimization algorithm.
\newpage

%

\end{document}